\documentclass[12pt]{article}
\usepackage{amsfonts}
\usepackage[dvips]{graphics}
\newtheorem{theorem}{Theorem}

\newtheorem{definition}{Definition}
\newtheorem{remark}{Remark}

\begin{document}
\title {A Liouville comparison principle for sub- and super-solutions of
the equation $w_t-\Delta_p (w) = |w|^{q-1}w$.}
\author{Vasilii V. Kurta.}
\maketitle
\begin{abstract}
\noindent We establish a Liouville comparison principle for entire
weak sub- and super-solutions of the equation $(\ast)$ $w_t-\Delta_p
(w) = |w|^{q-1}w$ in the half-space  ${\mathbb S}= {\mathbb
R}^1_+\times {\mathbb R}^n$, where $n\geq 1$, $q>0$ and  $ \Delta_p
(w):=\mbox{div}_x\left(|\nabla_x w|^{p-2}\nabla_x w\right)$,
$1<p\leq 2$. In our study we impose neither  restrictions on the
behaviour of entire weak  sub- and super-solutions on the
hyper-plane $t=0$,  nor any  growth conditions on their behaviour
and on that of any of their partial derivatives at infinity. We
prove that if $1<q\leq p-1+\frac pn$, and $u$ and $v$ are,
respectively, an entire weak super-solution  and an entire  weak
sub-solution of ($\ast$) in $\Bbb S$ which belong, only locally in
$\Bbb S$, to the corresponding Sobolev space and are such that
$u\geq v$, then $u\equiv v$. The result is sharp. As direct
corollaries we obtain known Fujita-type and Liouville-type theorems.
\end{abstract}
\thispagestyle{empty}
\section{Introduction and definitions.}

The purpose of this work is to obtain a  Liouville comparison
principle of  elliptic type for entire   weak sub- and
super-solutions of the equation
\begin{eqnarray}
w_t-\Delta_p (w) = |w|^{q-1}w
\end{eqnarray}
in the half-space  ${\mathbb S}= (0,+\infty)\times {\mathbb R}^n$,
where $n\geq 1$ is a natural number, $q>0$ is a  real number  and $
\Delta_p (w):=\sum_{i=1}^n \frac {d}{{dx}_i}A_i(\nabla w)$, with
$A_i(\xi)=|\xi|^{p-2}\xi_i$ for all $\xi=(\xi_1, \dots,
\xi_n)\in{\mathbb R}^n$  and $p>1$, defines the well-known
$p$-Laplacian operator. Under entire sub- and super-solutions of (1)
we understand sub- and super-solutions of (1) defined in the whole
half-space $\mathbb S$, and under Liouville results of elliptic type
for sub- and super-solutions  of the  parabolic equation (1)  in the
half-space $\Bbb S$ we understand Liouville-type results which, in
their formulations, have no restrictions on the behaviour of sub-
and super-solutions of (1) on the hyper-plane $t=0$. Also, we would
like to underline that we impose neither growth conditions on the
behaviour of sub- and super-solutions to (1) or on that of any of
their partial derivatives at infinity.

\begin{definition}
Let $n\geq 1$, $p>1$ and $q>0$. A  function $u=u(t,x)$ defined and
measurable in $\mathbb S$ is called  an entire   weak super-solution
of the equation (1) in $\mathbb S$ if it belongs to the function
space $ L_{q,\mathrm {loc}}(\mathbb S)$, with $u_t \in L_{1,\mathrm
{loc}}(\mathbb S)$ and $|\nabla_x u|^{p}\in L_{1, \mathrm
{loc}}(\mathbb S)$, and satisfies the integral inequality
\begin{eqnarray}
\int\limits_{\mathbb S}\left[u_t\varphi+\sum_{i=1}^n
|\nabla_x u|^{p-2}u_{x_i}\varphi_{x_i} -|u|^{q-1}u \varphi \right]dtdx
\geq 0
\end{eqnarray}
for every  non-negative function $\varphi \in  C^\infty (\mathbb S)$
with compact support in $\mathbb S$,  where $C^{\infty}({\mathbb
S})$ is the space of all functions defined and infinitely
differentiable in  ${\mathbb S}$.
\end{definition}
\begin{definition} A function $v=v(t,x)$ is an entire   weak sub-solution of (1)
if $u=-v$ is an entire   weak  super-solution of (1).
\end{definition}

\section{Results.}

\begin{theorem} Let $n\geq 1$,  $2\geq p> 1$
and  $1<q\leq p-1 +\frac {p}n$,  and let  $u$ be an entire   weak
super-solution and $v$ an entire   weak sub-solution of  (1) in
$\mathbb S$ such that $u\geq v$. Then $u = v$ in $ \mathbb S$.
\end{theorem}

The result in Theorem 1, which evidently has a comparison principle
character, we term a Liouville-type comparison principle, since, in
the particular cases when $u\equiv 0$ or $v\equiv 0$, it becomes a
Liouville-type theorem of elliptic type, respectively,  for entire
weak sub-solutions  or entire weak super-solutions of (1).

Since in Theorem 1 we impose no conditions on the behaviour of
entire weak  sub- or super-solutions of  the equation (1) on the
hyper-plane $t=0$, we can formulate,  as a direct corollary  of
Theorem 1, the following comparison principle, which in turn one can
term a Fujita comparison principle,  for entire weak sub- and
super-solutions of the Cauchy problem for the equation (1). It is
clear that in the particular cases when $u\equiv 0$ or $v\equiv 0$,
it becomes a Fujita-type theorem, respectively, for entire weak
sub-solutions  or entire weak super-solutions of the Cauchy problem
for the equation (1).

\begin{theorem} Let $n\geq 1$,  $2\geq p> 1$
and  $1<q\leq p-1 +\frac {p}n$, and let  $u$ be an entire   weak
super-solution and $v$ an entire   weak sub-solution of the Cauchy
problem, with possibly different initial data for $u$ and $v$, for
the equation (1) in $\mathbb S$ such that $u\geq v$. Then $u = v$ in
$\mathbb S$.
\end{theorem}
\begin{remark}
The  initial data for $u$ and $v$ in Theorem 2 may be different.
\end{remark}

Note that the results  in Theorems 1 and 2 are sharp and that the
hypotheses on the parameter $p$ in these theorems in fact force $p$
to be greater than $\frac{2n}{n+1}$. The sharpness of these results
for $q> p-1 +\frac {p}n\geq 1$ follows, for example, from the
existence of non-negative self-similar entire solutions to (1) in
$\Bbb S$, which was shown in [1]. Also, there one can find a
Fujita-type theorem on the non-existence  of non-negative entire
solutions of the Cauchy problem for (1), which was obtained as a
very interesting generalization of the famous blow-up result
established in [4], [5] and [9] to quasilinear parabolic equations.
For $0<q\leq1$, it is evident that the function $u(t,x)=e^t$  is a
positive entire classical super-solution of (1) in  $\Bbb S$.

We would also like to note that the results of the present work were
announced in [13] and that similar results for solutions of
semilinear parabolic inequalities were obtained in [7]. To prove the
results we use the $\alpha$-monotoni\-city property of the
$p$-Laplacian operator which was established in [11] and continue to
develop an approach in [7] and  [8], the elliptic analogue of which
was proposed [11]. That approach was subsequently used and developed
in the same framework by E. Mitidieri, S. Pokhozhaev and many
others, almost none of which cite the original research in [11].

For a survey of the literature on the asymptotic behaviour of and
blow-up results for solutions, sub- and super-solutions of the
Cauchy problem for nonlinear parabolic equations we refer to [2],
[3], [6], [14], [15] and [16].

\section {Proofs.}

In what follows, for  $q>1$ and   $2\geq p>1$,  let
\begin{eqnarray}\omega=\frac{p(q-1)}{q-p+1}\end{eqnarray}
and \begin{eqnarray}P (R)=\{(t,x)\in \mathbb S\!:\,
t^{2/\omega}+|x|^2<R^{2/\omega}\}\nonumber \end{eqnarray} for all
$R>0$. In this case it is clear that $0<\omega\leq 2$  and that  the
inequality
\begin{eqnarray}\hbox{volume of }P(R)\leq
{c}R^{\frac{n+\omega}\omega},\end{eqnarray} with $c$ some positive
constant which depends possibly only on $n$ and $\omega$, holds for
all $R>0$.

\noindent {\it Proof of Theorem 1.} Let $n\geq 1$,  $2\geq p> 1$ and
$1<q\leq p-1 +\frac {p}n$,  and let  $u$ be an entire   weak
super-solution and $v$ an entire  weak sub-solution of  (1) in
$\mathbb S$ such that $u\geq v$. By the well-known inequality
$$(|u|^{q-1}u-|v|^{q-1}v)(u-v)\geq 2^{1-q} |u-v|^{q+1}
$$
which holds for every $q\geq 1$ and all $u,v\in {\Bbb R}^1$  we
obtain from (2) the relation
\begin{eqnarray} \int\limits_{\mathbb S}\left [(u-v)_t\varphi
+\sum_{i=1}^n
{\varphi}_{x_i}(|\nabla_x u|^{p-2}u_{x_i}-
|\nabla_x v|^{p-2}v_{x_i}) \right]dtdx \geq \nonumber \\
2^{1-q}\int\limits_{\mathbb S}( u- v)^q  \varphi dtdx,
\end{eqnarray}
which holds for every non-negative function $\varphi\in
C^{\infty}(\Bbb S)$  with compact support in $\mathbb S$. Let
$\tau>0$ and $R>0$ be real numbers. Let $\eta\! :[0,+\infty)\to
[0,1]$ be a $C^\infty$-function   which has the non-negative
derivative $\eta'$ and equals 0 on the interval $[0, \tau]$ and 1 on
the interval $[2\tau, +\infty)$, and let $\zeta\! :[0,+\infty)\times
{\mathbb R}^n \to [0,1]$  be a $C^\infty$-function   which equals 1
on $\overline {P (R/2)}$ and 0 on $\{[0, +\infty) \times {\mathbb
R}^n\} \setminus \overline {P(R)}$. Let $\varphi (t,x)=
(w(t,x)+\varepsilon )^{-\nu}\zeta ^s(t,x)\eta^2(t)$,  where
$w(t,x)=u(t,x)-v(t,x)$, $\varepsilon>0$ and the positive constants
$s>1$ and $\nu\in (0, p-1)$ will be chosen below. Substituting the
function $\varphi$ in (5) and then integrating by parts   we arrive
at
\begin{eqnarray}
 -\frac s{1-\nu}
\int\limits_{P(R)}(w+\varepsilon)^{1-\nu}
\zeta_{t}\zeta^{s-1}\eta^2dtdx -\frac 2{1-\nu}
\int\limits_{P(R)}(w+\varepsilon)^{1-\nu} \zeta^s\eta' \eta
dtdx\nonumber
\\ -\nu \int\limits_{P(R)} \sum\limits_{i=1}^n
w_{x_i}
(|\nabla_x u|^{p-2}u_{x_i}-|\nabla_x v|^{p-2}v_{x_i})
(w+\varepsilon )^{-\nu -1}\zeta^s \eta^2 dtdx\nonumber\\
+ s\int\limits_{P(R)}
\sum_{i=1}^n
{\zeta}_{x_i}(|\nabla_x u|^{p-2}u_{x_i}-|\nabla_x v|^{p-2}v_{x_i})
(w+\varepsilon )^{-\nu}\zeta^{s-1} \eta^2 dtdx \nonumber
\\ \equiv I_1+I_2+I_3 +I_4 \geq 2^{1-q}\int\limits_{P(R)}w^q
(w+\varepsilon )^{-\nu}\zeta^s \eta^2 dtdx.
\end{eqnarray}
First, observing that $I_3$ is non-positive, we estimate $I_4$ in
terms of $I_3$ using the fact, which is a key point  in our proof,
that for $1<p\leq 2$ the $p$-Laplacian operator $\Delta_p$ satisfies
the $\alpha$-monotonicity condition (see, e.g, [12]) with
$\alpha=p$. This  in our case consists mostly of the fact that there
exists a positive constant $\mathcal K$ such that the coefficients
$A_i$, $i=1, \dots, n$, of the $p$-Laplacian operator satisfy the
inequality
\begin{eqnarray}
\left(\sum_{i=1}^n
\left(A_i\left(\xi^1\right)-A_i\left(\xi^2\right)\right)^2\right)
^{\alpha/2}
\leq {\mathcal K} \left(\sum_{i=1}^n\left(\xi_i^1-\xi_i^2\right)
\left(A_i\left(\xi^1\right)-A_i\left(\xi^2\right)\right) \right)
^{\alpha-1}\end{eqnarray} for all pairs $\xi^1, \xi^2 \in {\Bbb
R}^n$ and $\alpha=p$, provided $1<p\leq 2$. As a result, we have the
relation
\begin {eqnarray}
|I_4|\nonumber \\
\leq \int\limits_{P(R)} c_1|\nabla_x \zeta|\left(\sum_{i=1}^n
w_{x_i} (|\nabla_x u|^{p-2}u_{x_i}-|\nabla_x
v|^{p-2}v_{x_i})\right)^{\frac{p-1}p}
(w+\varepsilon)^{-\nu}\zeta^{s-1}\eta^2 dtdx. \end{eqnarray} Here we
use the symbols  $c_i$, $i=1,\dots, 8$, to denote constants
depending possibly   on $n$, $p$, $q$, $s$ or  $\nu$ but not on $R$,
$\varepsilon$ or $\tau$. Further, estimating   the integrand on the
right-hand side of (8)  by Young's inequality
\begin {eqnarray}
{\mathcal A}{\mathcal B}\leq \rho {\mathcal
A}^{\frac\beta{\beta-1}}+\rho^{1-\beta} {\mathcal B}^\beta
\end{eqnarray}
with $\rho=\frac {\nu}2, \beta=p,$
$${\mathcal A}=\left(\sum_{i=1}^n w_{x_i}  (|\nabla_x u|^{p-2}u_{x_i}-|\nabla_x v|^{p-2}v_{x_i})\right)^{\frac{p-1}p}
(w+\varepsilon)^{\frac{(1+\nu)(1-p)}p} \zeta^{\frac
{s(p-1)}p}\eta^\frac{2(p-1)}p$$ and
$${\mathcal B}=c_1|\nabla_x \zeta|
(w+\varepsilon)^{\frac {p-1-\nu}p}
\zeta^{\frac sp-1}\eta^\frac
2p, $$ we arrive at
\begin{eqnarray}
|I_4|\leq\frac {\nu}2  \int\limits_{P(R)}\sum_{i=1}^n w_{x_i}(|\nabla_x u|^{p-2}u_{x_i}-|\nabla_x v|^{p-2}v_{x_i})(w+\varepsilon)^{-\nu-1}\zeta^s\eta^2 dtdx
\nonumber \\
+\int\limits_{P(R)}c_2(w+\varepsilon)^{p-1-\nu}|\nabla_x\zeta|^p
\zeta^{s-p}\eta^2dtdx. \end{eqnarray} Now, observing that $I_2$ in
(6) is also non-positive, we obtain from (6) and (10) the relation
\begin{eqnarray} \int\limits_{P(R)}c_2
(w+\varepsilon)^{1-\nu}
|\zeta_{t}|\zeta^{s-1}\eta^2dtdx
+\int\limits_{P(R)}c_2(w+\varepsilon)^{p-1-\nu}|\nabla_x\zeta|^p\zeta^{s-p}\eta^2dtdx\nonumber\\
\geq \int\limits_{P(R)} w^q (w+\varepsilon)^{-\nu}\zeta^s \eta^2
dtdx\nonumber \\ +  \int\limits_{P(R)}\sum_{i=1}^n
w_{x_i}(|\nabla_x u|^{p-2}u_{x_i}-|\nabla_x
v|^{p-2}v_{x_i})(w+\varepsilon)^{-\nu-1} \zeta^s\eta^2 dtdx.
\end{eqnarray}
Estimating  both integrands on the left-hand side of (11) by Young's
inequality (9) with $\rho=\frac 12$, $\beta=\frac{q-\nu}{q-1},$
$${\mathcal A}=(w+\varepsilon)^{1-\nu}{\zeta}^{ \frac
{s(1-\nu)}{q-\nu}}\eta^{\frac{2(1-\nu)}{q-\nu}},$$ $${\mathcal
B}=c_2 |\zeta_t| {\zeta}^{{\frac
{s(q-1)}{q-\nu}}-1}\eta^{\frac{2(q-1)}{q-\nu}}$$  and $\rho=\frac
12$, $\beta=\frac{q-\nu}{q-p+1},$
$${\mathcal A}=(w+\varepsilon)^{p-1-\nu}{\zeta}^{\frac
{s(p-1-\nu)}{q-\nu}}\eta^{\frac{2(p-1-\nu)}{q-\nu}},$$
$${\mathcal B}= c_2|\nabla_x \zeta
|^p {\zeta}^{\frac {s(q-p+1)}{q-\nu}
-p}\eta^{\frac{2(q-p+1)}{q-\nu}},$$ respectively, we have the relation
\begin{eqnarray} \frac 12
\int\limits_{P(R)}(w+\varepsilon)^{q-\nu}\zeta^s \eta^2 dtdx
+c_3\int\limits_{P(R)}|\zeta_t|^{\frac {q-\nu}{q-1}}
\zeta^{s-\frac{q-\nu}{q-1}}\eta^2dtdx \nonumber \\ +\frac {1}2
\int\limits_{P(R)}(w+\varepsilon)^{q-\nu}\zeta^s \eta^2  dtdx + c_3
\int\limits_{P(R)} |\nabla_x \zeta|^{\frac
{p(q-\nu)}{q-p+1}}
\zeta^{s-{\frac{p(q-\nu)}{q-p+1}}}\eta^2dtdx \nonumber \\ \geq
\int\limits_{P(R)}w^q (w+\varepsilon)^{-\nu}\zeta^s \eta^2 dtdx
+
\nonumber \\   \int\limits_{P(R)}\sum_{i=1}^n w_{x_i}
(|\nabla_x u|^{p-2}u_{x_i}-|\nabla_x v|^{p-2}v_{x_i})(w+\varepsilon)^{-\nu-1}
\zeta^s\eta^2 dtdx.
\end{eqnarray}

Further, we  estimate the integral
$$
\int\limits_{P(R)}w^q\zeta^s \eta^2 dtdx
$$
by the inequality (12). To this end,  we substitute
$$\varphi (t,x) = \zeta^s(t,x) \eta^2(t)$$ in  (5) and after
integration by parts there we obtain
\begin{eqnarray} -s\int\limits_{P(R)}w\zeta_t\zeta^{s-1}\eta^2dtdx -2
\int\limits_{P(R)}w\zeta^s\eta' \eta dtdx \nonumber \\
+s\int\limits_{P(R)}\sum\limits_{i=1}^n
\zeta_{x_i}(|\nabla_x u|^{p-2}u_{x_i}-|\nabla_x v|^{p-2}v_{x_i}) \zeta^{s-1} \eta^2 dtdx
\nonumber \\
\geq 2^{1-q}\int\limits_{P(R)}w^q\zeta^s \eta^2 dtdx.
\end{eqnarray}
Since  the second term  on the left-hand side of (13) is
non-positive we have
\begin{eqnarray}
s\int\limits_{P(R)}w|\zeta_t|\zeta^{s-1}\eta^2dtdx \nonumber \\
+s\int\limits_{P(R)}\sum\limits_{i=1}^n \zeta_{x_i} (|\nabla_x
u|^{p-2}u_{x_i}-|\nabla_x v|^{p-2}v_{x_i}) \zeta^{s-1} \eta^2 dtdx
\geq 2^{1-q} \int\limits_{P(R)}w^q\zeta^s \eta^2 dtdx.
\end{eqnarray} Now, estimating the first integral on the left-hand
side of (14) by H\"{o}lder's  inequality, we arrive at
\begin{eqnarray} s\left(\int\limits_{P(R)\setminus P(R/2)}w^q\zeta^s
\eta^2 dtdx\right)^{\frac 1q}
\left(\int\limits_{P(R)}|\zeta_t|^{\frac q
{q-1}} \zeta^{s-\frac q{q-1}} \eta^2 dtdx\right)^{\frac {q-1}q} \nonumber \\
+s\int\limits_{P(R)}\sum\limits_{i=1}^n \zeta_{x_i}(|\nabla_x
u|^{p-2}u_{x_i}-|\nabla_x v|^{p-2}v_{x_i}) \zeta^{s-1} \eta^2 dtdx
\nonumber
\\ \geq 2^{1-q}\int\limits_{P(R)}w^q\zeta^s \eta^2 dtdx.
\end{eqnarray}
On the other hand, by (7) we have
\begin{eqnarray}
\int\limits_{P(R)}\sum\limits_{i=1}^n\zeta_{x_i}(|\nabla_x u|^{p-2}u_{x_i}-|\nabla_x v|^{p-2}v_{x_i})
\zeta^{s-1} \eta^2 dtdx\nonumber \\ \leq c_4 \int\limits_{P(R)} |\nabla_x\zeta|
\left(w_{x_i}(|\nabla_x u|^{p-2}u_{x_i}-|\nabla_x v|^{p-2}v_{x_i}) \right)^{\frac{p-1}p}
\zeta^{s-1}\eta^2 dtdx.
\end{eqnarray}
Estimating  the right-hand side of (16) by H\"{o}lder's inequality
we arrive at the relation
\begin{eqnarray}
\int\limits_{P(R)}\sum\limits_{i=1}^n\zeta_{x_i}(|\nabla_x u|^{p-2}u_{x_i}-|\nabla_x v|^{p-2}v_{x_i})
\zeta^{s-1} \eta^2 dtdx\nonumber \\ \leq c_4 \left(\int\limits_{P(R)}(w+\varepsilon )^
{(1+\nu)(p-1)}|\nabla_x \zeta|^p \zeta^{s-p}\eta^2
dtdx\right) ^{1/p}\nonumber \\ \times
\left(\int\limits_{P(R)}\sum_{i=1}^n w_{x_i}(|\nabla_x u|^{p-2}u_{x_i}-|\nabla_x v|^{p-2}v_{x_i})
(w+\varepsilon)^{-\nu-1}\zeta^s\eta^2
dtdx\right)^{\frac{p-1}p}
\end{eqnarray}
which holds for every $\varepsilon >0$ and $p-1> \nu>0$. Further,
for any $d>1$ we have
\begin{eqnarray}\int\limits_{P(R)} (w+\varepsilon
)^{(1+\nu)(p-1)}
|\nabla_x \zeta|^p \zeta^{s-p}\eta^2 dtdx\nonumber \\ \leq
\left(\int\limits_{P(R)\setminus  P(R/2)}
(w+\varepsilon)^{d(1+\nu)(p-1)}\zeta^s \eta^2 dtdx
\right)^{\frac 1 d}\nonumber \\
\times
\left (\int\limits_{P(R)}|\nabla_x \zeta|^{\frac {p d}{d-1}}
\zeta^{s-\frac{p d}{d-1}}\eta^2 dtdx\right)^{\frac {d-1} d}.
\end{eqnarray}
Now, we choose for every $q>1$ and a sufficiently small $\nu$ from
the interval $(0, p-1)$ the parameter $d$ such that
$d(1+\nu)(p-1)=q$. Then (17) and (18) yield
\begin{eqnarray}
\int\limits_{P(R)}\sum\limits_{i=1}^n\zeta_{x_i}(|\nabla_x u|^{p-2}u_{x_i}-|\nabla_x v|^{p-2}v_{x_i}) \zeta^{s-1} \eta^2 dtdx\nonumber \\
\leq c_4\left(\int\limits_{P(R)\setminus  P(R/2)}
(w+\varepsilon)^q\zeta^s \eta^2 dtdx \right)^{\frac 1 {p  d}}
\left (\int\limits_{P(R)}|\nabla_x \zeta|^{\frac {p d}{d-1}}
\zeta^{s-\frac{p d}{d-1}}\eta^2 dtdx\right)^{\frac {d-1}
{p  d}}\nonumber \\ \times \left(\int\limits_{P(R)}\sum_{i=1}^n
w_{x_i}(|\nabla_x u|^{p-2}u_{x_i}-|\nabla_x v|^{p-2}v_{x_i})(w+\varepsilon)^{-\nu-1}\zeta^s\eta^2
dtdx\right)^{\frac{p-1} p}.
\end{eqnarray}
Estimating the last term on the right-hand side of (19) by (12), we
have   \begin{eqnarray}
\int\limits_{P(R)}\sum\limits_{i=1}^n\zeta_{x_i}(|\nabla_x
u|^{p-2}u_{x_i}-|\nabla_x v|^{p-2}v_{x_i}) \zeta^{s-1} \eta^2
dtdx\nonumber \\
\leq c_4\left(\int\limits_{P(R)\setminus P(R/2)}
(w+\varepsilon)^q\zeta^s \eta^2 dtdx \right)^{\frac 1{p d}} \left
(\int\limits_{P(R)}|\nabla_x \zeta|^{\frac {p d}{d-1}}
\zeta^{s-\frac{p d}{d-1}}\eta^2 dtdx\right)^{\frac {d-1} {p
d}}\nonumber
\\
\times\left(
\int\limits_{P(R)}(w+\varepsilon)^{q-\nu}\zeta^s \eta^2 dtdx\right. -\int\limits_{P(R)}w^q(w+\varepsilon)^{-\nu}\zeta^s\eta^2 dtdx \nonumber \\
+c_3 \int\limits_{P(R)}|\zeta_t|^{\frac{q-\nu}{q-1}} \zeta^{s-\frac
{q-\nu}{q-1}}\eta^2dtdx  +c_3 \left. \int\limits_{P(R)} |\nabla_x
\zeta|^{\frac{p(q-\nu)}{q-p+1}}
\zeta^{s-\frac{p(q-\nu)}{q-p+1}}\eta^2dtdx \right)^\frac{p-1}p.
\end{eqnarray}
In (20), passing to the limit as $\varepsilon\to 0$ as justified  by
Lebesgue's theorem (see, e.g., [10, p. 303]) we obtain
\begin{eqnarray}
\int\limits_{P(R)}\sum\limits_{i=1}^n\zeta_{x_i}(|\nabla_x u|^{p-2}u_{x_i}-|\nabla_x v|^{p-2}v_{x_i})
\zeta^{s-1} \eta^2 dtdx\nonumber \\ \leq
c_5\left(\int\limits_{P(R)\setminus P(R/2)}w^q\zeta^s \eta^2 dtdx
\right)^{\frac 1{p d}} \left (\int\limits_{P(R)}|\nabla_x
\zeta|^{\frac {p d}{d-1}}\zeta^{s-\frac{p d}{d-1}}\eta^2
dtdx\right)^{\frac {d-1} {p d}} \nonumber \\ \times
\left(\int\limits_{P(R)}|\zeta_t|^{\frac {q-\nu}{q-1}}\zeta^{s-\frac
{q-\nu}{q-1}}\eta^2 dtdx + \int\limits_{P(R)} |\nabla_x
\zeta|^{\frac
{p(q-\nu)}{q-p+1}}\zeta^{s-\frac{p(q-\nu)}{q-p+1}}dtdx\right)^{\frac{p-1}p}.
\end{eqnarray}
Further, (15) and (21)  yield
\begin{eqnarray} \int\limits_{P(R)}w^q\zeta^s \eta^2 dtdx\leq
c_6\left(\int\limits_{P(R)\setminus P(R/2)}w^q\zeta^s \eta^2
dtdx\right)^{\frac 1q} \left(\int\limits_{P(R)}|\zeta_t|^{\frac q
{q-1}} \zeta^{s-\frac q{q-1}}\eta^2 dtdx\right)^{\frac {q-1}q}\nonumber \\
+c_6\left(\int\limits_{P(R)\setminus P(R/2)} w^q \zeta^s \eta^2 dtdx
\right)^{\frac 1 {p d}} \left(\int\limits_{P(R)}|\nabla_x
\zeta|^{\frac{p  d}{d-1}}\zeta^{s-\frac{p d}{d-1}}\eta^2  dtdx\right)^{\frac {d-1} {p d}} \nonumber \\
\times \left(\int\limits_{P(R)}|\zeta_t|^{\frac {q-\nu}{q-1}}\zeta^{s-\frac {q-\nu}{q-1}}\eta^2 dtdx
+\int\limits_{P(R)}|\nabla_x
\zeta|^{\frac{p(q-\nu)}{q-p+1}}
\zeta^{s-\frac{p(q-\nu)}{q-p+1}}\eta^2 dtdx\right)^{\frac
{p-1}p}.
\end{eqnarray}
Now, for arbitrary $(t,x)\in \mathbb S$ and $R>0$, we choose in (22)
the function $\zeta=\zeta(t, x)$ in the form
\begin{equation}
\zeta
(t,x)=
\psi\left(\frac {t^{2/\omega} +|x|^2}{R^{2/\omega}}\right),
\end{equation}
where $0<\omega\leq 2$ is given by (3)  and   $\psi: [0, \infty)\to
[0,1]$ is a $C^\infty$-function   which equals 1 on $[0, 2^{-\frac
2\omega}]$ and 0 on $[1, \infty)$ and is such that the inequalities
\begin{eqnarray}
|\zeta_t|\leq {c_7}R^{-1}\qquad \hbox{and}\qquad  |\nabla_x\zeta|\leq
{c_7}R^{-\frac {1}\omega}\end {eqnarray}
hold.
Note that
it is always possible to find such a function $\zeta$. Indeed, this
can be easily verified by direct calculation of the corresponding
derivatives  of the function $\zeta$  given
by (23).
Also, choosing in  (22) the parameter $s$ sufficiently large, we have from (22) by  (4) and (24) the relation
\begin{eqnarray} \int\limits_{P(R)}w^q\zeta^s \eta^2 dtdx \leq
c_8\left(R^{\frac {n+\omega}\omega-\frac q
{q-1}}\right)^{\frac{q-1}q} \left(\int\limits_{P(R)\setminus
P(R/2)}w^q\zeta^s \eta^2 dtdx\right)^{\frac  1q}\nonumber \\ +
c_8\left(R^{\frac {n+\omega}\omega-\frac {p d} {\omega(d-1)}}\right)^{\frac{d-1}{p d}} \left(R^{\frac
{n+\omega}\omega-\frac {q-\nu} {q-1}}+  R^{\frac
{n+\omega}\omega-\frac {p(q-\nu)} {\omega(q-p+1)}} \right)^{\frac
{p-1}p} \left(\int\limits_{P(R)\setminus P(R/2)} u^q \zeta^s
\eta^2 dtdx \right)^{\frac 1 {p d}},
\nonumber \end{eqnarray}
which in turn by (3)
implies
\begin{eqnarray} \int\limits_{P(R)}w^q\zeta^s \eta^2 dtdx \leq
c_8\left(R^{\frac {n+\omega}\omega-\frac q
{q-1}}\right)^{\frac{q-1}q} \left(\int\limits_{P(R)\setminus
P(R/2)}w^q\zeta^s \eta^2 dtdx\right)^{\frac  1q}\nonumber \\ +
c_8\left(R^{\frac {n+\omega}\omega-\frac {p d} {\omega(d-1)}}\right)^{\frac{d-1}{p d}} \left(R^{\frac
{n+\omega}\omega-\frac {q-\nu} {q-1}} \right)^{\frac
{p-1}p} \left(\int\limits_{P(R)\setminus P(R/2)} u^q \zeta^s
\eta^2 dtdx \right)^{\frac 1 {p d}}.
\end{eqnarray}
Making simple calculation in (25) we arrive at
\begin{eqnarray}
\int\limits_{P(R)}w^q\zeta^s \eta^2
dtdx \leq c_8R^{\frac n{p q}\left[q-p+1 -\frac p
n\right]} \left(\int\limits_{P(R)\setminus P(R/2)} w^q \zeta^s \eta^2
dtdx \right)^{\frac 1q}\nonumber \\ +c_8R^{\frac {n(p
q-p+1-\nu(p-1)} {p^2 q(q-1)}\left[q-p+1
- \frac p n\right]} \left(\int\limits_{P(R)\setminus P(R/2)} w^q
\zeta^s \eta^2 dtdx \right)^{\frac 1 {p d}}.
\end{eqnarray}
Further, since for $q>1$, $2\geq p>1$ and  $p-1>\nu>0$ the
quantities
$$
\frac n{p q}\qquad \hbox{and}\qquad \frac {n(p
q-p+1-\nu(p -1))}{p^2 q(q-1)}
$$
are positive, we obtain from (26) for
$1 <q<p-1+\frac p n
$
the relation  \begin{equation}\int_{\mathbb S} u^q \eta^2dtdx=0.\end{equation}
Also, for  $q = p-1 + \frac p n$ we deduce  from (26) that
$$\int_{\mathbb S} u^q \eta^2dtdx<\infty.$$ The latter  yields
the relation
\begin{eqnarray}
\int\limits_{P(R_k)\setminus P({R_k}/2)} w^q \eta^2 dtdx\to 0
\end{eqnarray}
which holds for every  sequence $R_k\to \infty$. On the other hand, the inequality  \begin{eqnarray} \int\limits_{P(R/2)}w^q \eta^2 dtdx
\leq c_8R^{\frac n{p q}\left[q-p+1- \frac p n\right]}
\left(\int\limits_{P(R)\setminus P(R/2)} w^q \eta^2 dtdx
\right)^{\frac 1q}\nonumber \\ +c_8R^{\frac {n(p
q-p+1-\nu(p -1))} {p^2 q(q-1)}\left[q- p+1-\frac
p n\right]} \left(\int\limits_{P(R)\setminus P(R/2)} w^q \eta^2
dtdx \right)^{\frac 1 {p d}}
\end{eqnarray}
follows easily from (26).  In turn,  (28) and (29) imply for
$q=p-1+\frac p n$ that the relation
\begin{eqnarray} \int\limits_{P(R_k)}u^q \eta^2dtdx \to 0\nonumber
\end{eqnarray}
holds for every  sequence $R_k\to \infty$. The latter implies that
(27) holds for every $q$ satisfying \begin{equation} 1 <q\leq
p-1+\frac p n.\end{equation} In (27), by letting  the parameter
$\tau$ in the definition of the function $\eta$ tend to zero, we
obtain that $u(t,x)=v(t,x)$ a.e. in $\mathbb S$ for every $q$ which
satisfies (30). $\Box$

\vskip 10pt \noindent {\bf Acknowledgments.}

\vskip 10pt \noindent This research was financially supported by the
Alexander von Humboldt Foundation (AvH). The author is very grateful
to AvH  for the  opportunity to visit the Mathematical Institute of
K\"oln University and to Professor B. Kawohl for his cordial
hospitality during this visit.

\vspace{10mm}

\noindent \textbf{Author's address:}

\vspace{10 mm}

\noindent Vasilii V. Kurta

\noindent American Mathematical Society

\noindent Mathematical Reviews

\noindent 416 Fourth Street, P.O. Box 8604

\noindent Ann Arbor, Michigan 48107-8604, USA

\noindent \textbf {e-mail:} vkurta@umich.edu, vvk@ams.org
\end{document}